\documentclass[12pt]{amsart}

\usepackage{amsfonts,amsmath,amssymb,amsthm}
\usepackage{fancyhdr}
\usepackage[all]{xypic}

\newcommand{\Cc}{\mathbb{C}}  
\newcommand{\Pp}{\mathbb{P}}
\newcommand{\Rr}{\mathbb{R}}

\newcommand{\defi}[1]{\emph{#1}}

\renewcommand{\epsilon}{\varepsilon}
\renewcommand{\le}{\leqslant}
\renewcommand{\ge}{\geqslant}
\renewcommand {\leq}{\leqslant}
\renewcommand {\geq}{\geqslant}

\newcommand{\grad}{\mathop{\mathrm{grad}}\nolimits}

\newcommand{\Baff}{{\mathcal{B}_\mathit{\!aff}}}
\newcommand{\Binf}{{\mathcal{B}_\infty}}
\newcommand{\B}{{\mathcal{B}}}

\newcommand{\Aut}{\mathop{\mathrm{Aut}}\nolimits}
\newcommand{\Disc}{\mathop{\mathrm{Disc}}\nolimits}

\newcommand{\id}{\mathop{\mathrm{id}}\nolimits}

\newcommand{\Diff}{\mathop{\mathrm{Diff}}\nolimits}
\newcommand{\ine}{\mathrm{in}}
\newcommand{\out}{\mathrm{out}}
\newcommand{\dic}{\mathrm{dic}}
{\theoremstyle{plain}
\newtheorem*{theorem*}{Theorem}  
\newtheorem{theorem}{Theorem}    

\newtheorem{lemma}[theorem]{Lemma}        

}

{\theoremstyle{remark}
\newtheorem*{remark*}{Remark}  
\newtheorem{remark}{Remark}  
\newtheorem{example}{Example}
}

\begin{document}
{\large\title{\textsf{Invariance of Milnor numbers \\ and  topology of complex polynomials}}}
\author{\textbf{\textsf{Arnaud Bodin}}}
\begin{abstract}
We give a global version of L\^e-Ramanujam $\mu$-constant theorem for
polynomials. Let $(f_t)$, ${t\in[0,1]}$, be a family of polynomials of $n$
complex variables with isolated singularities, whose coefficients are
polynomials in  $t$. We consider the case where some numerical invariants
are constant (the affine Milnor number $\mu(t)$, the  Milnor number at
infinity $\lambda(t)$, the number of critical values, the number of
affine critical values, the number of critical values at infinity). Let
$n=2$, we also suppose the degree of the $f_t$ is a constant, then the
polynomials $f_0$ and $f_1$ are topologically equivalent. For $n>3$ we
suppose that critical values at infinity depend continuously on $t$, then
we prove that the geometric monodromy representations of the $f_t$, are
all equivalent.
\end{abstract}

\maketitle


\section{Introduction}

Let $f : \Cc^n \longrightarrow \Cc$ be a polynomial map, $n\ge2$. By a
result of Thom \cite{T} there is a minimal
 \defi{set of critical values} $\B$ of point of $\Cc$ such that
$f : f^{-1}(\Cc\setminus \B) \longrightarrow \Cc\setminus \B$
is a fibration.

\subsection{Affine singularities}
We suppose that \defi{affine singularities are isolated} \emph{i.e.} that the
set $\{ x \in \Cc^n \ | \ \grad_f x = 0\}$ is a finite set.
Let $\mu_c$ be the sum of the local Milnor numbers at the points of $f^{-1}( c )$.
Let
$$\Baff = \big\{ c \ | \ \mu_c > 0\big\} \quad \text { and }\quad \mu = \sum_{c\in\Cc} \mu_c$$
be the \defi{affine critical values} and the \defi{affine Milnor number}.

\subsection{Singularities at infinity}

See \cite{Br}. Let $d$ be the degree of $f : \Cc^n \longrightarrow \Cc$,
let $f = f^d+f^{d-1}+\cdots+f^0$ where $f^j$ is homogeneous of degree $j$.
Let $\bar f (x,x_0)$ (with $x = (x_1,\ldots,x_n)$) be the homogenization
of $f$ with the new variable $x_0$: $\bar f(x,x_0) =
f^d(x)+f^{d-1}(x)x_0+\ldots +f^0(x)x_0^d$. Let $$X = \left\lbrace
((x:x_0),c)\in \Pp^n\times \Cc \ | \ \bar f(x,x_0)-cx_0^d=0 \right\rbrace.
$$ Let $\mathcal{H}_\infty$ be the hyperplane at infinity of $\Pp^n$
defined by $(x_0=0)$. The singular locus of $X$ has the form $\Sigma
\times \Cc$ where $$\Sigma = \left\lbrace (x:0) \ | \ \frac{\partial
f^d}{\partial x_1}=\cdots =\frac{\partial f^d}{\partial x_n}= f^{d-1} =0
\right\rbrace \subset \mathcal{H}_\infty.$$ We suppose that $f$ has
\defi{isolated singularities at infinity} that is to say that $\Sigma$ is
finite. This is always true for $n=2$. For a point $(x:0) \in
\mathcal{H}_\infty$, assume, for example, that $x =
(x_1,\ldots,x_{n-1},1)$ and set $\check{x} = (x_1,\ldots,x_{n-1})$ and
$$F_c(\check{x},x_0) = \bar f(x_1,\ldots,x_{n-1},1)-cx_0^d.$$ Let
$\mu_{\check{x}}(F_c)$ be the local Milnor number of $F_c$ at the point
$(\check{x},0)$. If $(x:0)\in \Sigma$ then $\mu_{\check{x}}(F_c)>0$. For
a generic $s$,  $\mu_{\check{x}}(F_s) =  \nu_{\check{x}}$, and for
finitely many $c$, $\mu_{\check{x}}(F_c) > \nu_{\check{x}}$. We set
$\lambda_{c,\check{x}} = \mu_{\check{x}}(F_c) - \nu_{\check{x}}$,
$\lambda_c = \sum_{(x:0)\in\Sigma} \lambda_{c,\check{x}}$. Let $$\Binf =
\big\{ c\in\Cc \ | \  \lambda_c > 0\big\} \quad \text { and } \quad
\lambda = \sum_{c\in\Cc} \lambda_c$$ be the \defi{critical values at
infinity} and the \defi{Milnor number at infinity}. We can now describe
the set of critical values $\B$ as follows (see \cite{HL} and \cite{Pa}):
$$\B = \Baff \cup \Binf.$$ Moreover by \cite{HL} and \cite{ST} for
$s\notin \B$, $f^{-1}(s)$ has the homotopy type of a wedge of
$\lambda+\mu$ spheres of real dimension $n-1$.

\subsection{Statement of the results}

\begin{theorem}
\label{th:A} Let $(f_t)_{t\in[0,1]}$ be a family of complex polynomials
from $\Cc^n$ to $\Cc$ whose coefficients are polynomials in  $t$. We
suppose that affine singularities and singularities at infinity are
isolated. Let suppose that the integers $\mu(t)$, $\lambda(t)$,
$\#\B(t)$, $\#\Baff(t)$,  $\#\Binf(t)$ do not depend on $t \in[0,1]$.
Moreover let us suppose that critical values at infinity $\Binf(t)$
depend continuously on $t$. Then the fibrations  $f_0 :
f_0^{-1}(\Cc\setminus\B(0)) \longrightarrow \Cc\setminus\B(0)$ and $f_1 :
f_1^{-1}(\Cc\setminus\B(1)) \longrightarrow \Cc\setminus\B(1) $ are fiber
homotopy equivalent, and for $n\not=3$ are differentiably isomorphic.
\end{theorem}

\begin{remark}
As a consequence for $n\not= 3$ and $*\notin  \B(0)\cup  \B(1)$ the
monodromy representations $$\pi_1(\Cc \setminus \B(0),*) \longrightarrow
\Diff (f_0^{-1}(*)) \text{ and }$$ $$\pi_1(\Cc \setminus \B(1),*)
\longrightarrow \Diff (f_1^{-1}(*))$$ are equivalent (where $\Diff
(f_t^{-1}(*))$ denotes the diffeomorphisms of $f_t^{-1}(*)$ modulo
diffeomorphisms isotopic to identity).
\end{remark}

\begin{remark}
The restriction $n\not=3$, as in \cite{LR}, is due to the use of the $h$-cobordism theorem.
\end{remark}

\begin{remark}
This result extends a theorem of H\`a~H.V and Pham~T.S. \cite{HP} which
deals only with monodromy at infinity (which correspond to a loop around
the whole set $\B(t)$) for $n=2$. For $n>3$ the invariance of monodromy
at infinity is stated by M.~Tib\u{a}r in \cite{Ti}. The proof is based on
the articles of H\`a~H.V.-Pham~T.S. \cite{HP} and of
L\^e~D.T.-C.P.~Ramanujam \cite{LR}.
\end{remark}

\begin{lemma}
Under the hypotheses of the previous theorem (except the
hypothesis of continuity of the critical values), and one of the
following conditions:
\begin{itemize}
\item $n=2$, and $\deg f_t$ does not depend on t;
\item $\deg f_t$, and $\Sigma(t)$ do not depend on t, and for all
$(x:0)\in \Sigma(t)$, $\nu_{\check x}(t)$ is independent of $t$;
\end{itemize}
 we have that $\Binf(t)$ depends continuously on t,
 \emph{i.e.} if $c(\tau) \in \Binf(\tau)$ then for all $t$ near $\tau$ there
exists $c(t)$ near $c(\tau)$ such that $c(t) \in \Binf(t)$.
\end{lemma}

Under the hypothesis that there is no singularity at infinity we can prove
the stronger result:
\begin{theorem}
\label{th:B} Let $(f_t)_{t\in[0,1]}$ be a family of complex polynomials
whose coefficients are polynomials in $t$. Suppose that $\mu(t)$,
$\#\Baff(t)$ do not depend on $t \in[0,1]$. Moreover suppose that
$n\not=3$ and for all $t\in [0,1]$ we have $\Binf(t) = \varnothing$. Then
the polynomials $f_0$ and $f_1$ are topologically equivalent that is to
say there exists homeomorphisms $\Phi$ and $\Psi$ such that $$ \xymatrix{
{}\Cc^n \ar[r]^-\Phi\ar[d]_-{f_0}  & \Cc^n \ar[d]^-{f_1} \\
  \Cc \ar[r]_-\Psi        & \Cc.
}
$$
\end{theorem}

For the proof we glue the former study with the version
of the $\mu$-constant theorem of L\^e~D.T. and C.P.~Ramanujam
stated by J.G.~Timourian \cite{Tim}: a $\mu$-constant deformation
of germs of isolated hypersurface singularity is a product family.

For polynomials in two variables we can prove the following theorem which
is  a global version of L\^e-Ramanujam-Timourian theorem:
\begin{theorem}
\label{th:C} Let $n=2$. Let $(f_t)_{t\in[0,1]}$ be a family of complex
polynomials whose coefficients are polynomials in $t$. Suppose that the
integers $\mu(t)$, $\lambda(t)$, $\#\B(t)$, $\#\Baff(t)$,  $\#\Binf(t)$,
$\deg f_t$ do not depend on $t \in[0,1]$.  Then the polynomials $f_0$ and
$f_1$ are topologically equivalent.
\end{theorem}
It uses a result of L.~Fourrier \cite{Fo} that give a necessary and sufficient condition
for polynomials to be topologically equivalent outside sufficiently
large compact sets of $\Cc^2$.

\bigskip
 This work was initiated by an advice of L\^e~D.T. concerning the
 article \cite{Bo}:``It is easier to find conditions for polynomials
 to be equivalent than find all polynomials that respect a given
 condition."

\bigskip

We will denote $B_R = \big\lbrace x\in\Cc^n \ | \ \|x\| \le R \big\rbrace$,
$S_R = \partial B_R =  \big\lbrace x\in\Cc^n \ | \ \|x\| = R \big\rbrace$
and $D_r(c) = \big\lbrace s\in\Cc \ | \ \|s-c\| \le r \big\rbrace$.

\section{Fibrations}

In this paragraph we give some properties for a complex polynomial
$f : \Cc^n \longrightarrow \Cc$.
The two first lemmas are consequences of transversality properties.
There are direct generalizations of lemmas of \cite{HP}.
Let $f$ be a polynomial of $n$ complex variables with isolated affine singularities
and with isolated singularities at infinity.
For each fiber $f^{-1}(c)$ there is a finite number of real numbers
$R>0$ such that $f^{-1}(c)$ has non-transversal intersection with the sphere
$S_R$. So for a sufficiently large number $R(c)$ the intersection
$f^{-1}(c)$ with $S_R$ is transversal for all $R \ge R(c)$.
Let $R_1$ be the maximum of the $R(c)$ with $c\in \B$.
We choose a small $\epsilon$,  $0<\epsilon \ll 1$ such that for all values $c$
in the bifurcation set $\B$ of $f$ and for all $s \in D_{\epsilon}(c)$
the intersection $f^{-1}(s) \cap S_{R_1}$ is transversal, this is possible by
continuity of the transversality. Let choose $r>0$ such that
$\B$ is contained in the interior of $D_r(0)$.
We denote
$$K = D_r(0) \setminus \bigcup_{c\in\B} \mathring D_{\epsilon}(c).$$

\begin{lemma}
\label{lem:1}
There exists $R_0 \gg 1$ such that for all $R \ge R_0$ and for all $s$ in $K$,
$f^{-1}(s)$ intersects $S_R$ transversally.
\end{lemma}

\begin{proof}
We have to adapt the beginning of the proof of \cite{HP}. If the
assertion is false then we have a sequence $(x_k)$ of points of $\Cc^n$
such that $f(x_k) \in K$ and $\| x_k \| \rightarrow +\infty$ as $k
\rightarrow +\infty$ and such that there exists complex numbers
$\lambda_k$ with $\grad_f x_k = \lambda_k x_k$, where the gradient is
Milnor gradient: $\grad_f = \left(\overline{\frac{\partial f}{\partial
x_1}},\ldots, \overline{\frac{\partial f}{\partial x_n}}\right)$. Since
$K$ is a compact set we can suppose (after extracting a sub-sequence, if
necessary) that $f(x_k) \rightarrow c\in K$ as $k \rightarrow +\infty$.
Then by the Curve Selection Lemma of \cite{NZ} there exists a real
analytic curve $x : ]0,\epsilon[ \longrightarrow \Cc^n$
such that $x(\tau) = a\tau^\beta + a_1\tau^{\beta+1}+\cdots$ with $\beta
< 0$, $a\in \Rr^{2n} \setminus \{ 0 \}$ and $\grad_f x(\tau) =
\lambda(\tau) x(\tau)$. Then $f(x(\tau)) = c + c_1\tau^\rho+\cdots$ with
$\rho > 0$. Then we can redo the calculus of \cite{HP}:
$$\frac{df(x(\tau))}{d\tau}= \big\langle \frac{dx}{d\tau},\grad_f x(\tau)
\big\rangle = \bar \lambda(\tau)  \big\langle \frac{dx}{d\tau}, x(\tau)
\big\rangle$$ it implies $$ |\lambda(\tau)| \leq 2 \frac
{\big|\frac{df(x(\tau))}{d\tau} \big|} {\frac{d\|x(\tau)\|^2}{d\tau}}.$$
As $\| x(\tau) \| = b_1\tau^\beta+\cdots$ with $b_1 \in \Rr_+^*$ and
$\beta < 0$ we have $|\lambda(\tau)| \leq \gamma
\frac{\tau^{\rho-1}}{\tau^{2\beta-1}}= \gamma \tau^{\rho-2\beta}$ where
$\gamma$ is a constant. We end the proof be using the characterization of
critical value at infinity in \cite{Pa}: $$ \|x(\tau)\|^{1-1/N}\| \grad_f
x(\tau) \| =  \|x(\tau)\|^{1-1/N}  |\lambda(\tau)|\,  \|x(\tau)\| \leq
\gamma \tau^{\rho-\beta/N}$$ As $\rho>0$ and $\beta<0$, for all $N>0$ we
have that $ \|x(\tau)\|^{1-1/N}\| \grad_f x(\tau) \| \rightarrow 0$ as
$\tau \rightarrow 0$. By \cite{Pa} it implies that the value $c$ (the
limit of $f(x(\tau))$ as $\tau \rightarrow 0$) is in $\Binf$. But as
$c\in K$ it is impossible.
\end{proof}

This first lemma enables us to get the following result:  because
of the transversality we can find a vector field tangent to the fibers
of $f$ and pointing out the spheres $S_R$. Integration of such a vector
field gives the next lemma.
\begin{lemma}
\label{lem:2}
The fibrations $f :  f^{-1}(K)\cap  \mathring B_{R_0} \longrightarrow K$ and
$f : f^{-1}(K) \longrightarrow K$ are differentiably isomorphic.
\end{lemma}

We will also need the following fact:
\begin{lemma}
\label{lem:0}
The fibrations $f :  f^{-1}(K) \longrightarrow K$ and
$f : f^{-1}(\Cc\setminus\B) \longrightarrow \Cc\setminus\B$ are differentiably isomorphic.
\end{lemma}

The following lemma is adapted from \cite{LR}.
For completeness we give the proof.
\begin{lemma}
\label{lem:3}
Let $R,R'$ with  $R\ge R'$ be real numbers such that the intersections
$f^{-1}(K)\cap  S_{R}$ and  $f^{-1}(K)\cap  S_{R'}$ are transversal.
Let us suppose that $f :  f^{-1}(K)\cap B_{R'}  \longrightarrow K$
and $f :   f^{-1}(K)\cap  B_{R} \longrightarrow K$
are fibrations with fibers homotopic to a wedge of $\nu$ $(n-1)$-dimensional spheres.
Then the fibrations are fiber homotopy equivalent.
And for $n \not= 3$
 the fibrations are differentiably equivalent.
\end{lemma}

\begin{proof}
The first part is a consequence of a result of A.~Dold \cite[th.~6.3]{Do}.
The first fibration is contained in the second. By the result of Dold we only have to prove
that if $* \in \partial D_r$ then the inclusion of
$F' = f^{-1}(*) \cap B_{R'}$ in $F = f^{-1}(*) \cap B_{R}$ is an homotopy equivalence.
To see this we choose a generic $x_0$ in $\Cc^n$ such that the real function
$x \mapsto \| x -x_0 \|$ has non-degenerate critical points of index less than $n$
(see \cite[\S7]{Mi}).
Then $F$ is obtained from
$F'$ by attaching cells of index less than $n$.
For $n=2$ the fibers are homotopic to a wedge of $\nu$ circles, then the inclusion
of $F'$ in $F$ is an homotopy equivalence.
For $n>2$ the fibers $F,F'$ are simply connected and the morphism
$H_i(F') \longrightarrow H_i(F)$ induced by inclusion is an isomorphism.
For $i\not=n-1$ this is obvious since $F$ and $F'$ have the homotopy type of a wedge of
$(n-1)$-dimensional spheres, and for $i=n-1$ the exact sequence of the pair $(F,F')$
is
$$ H_{n}(F,F') \longrightarrow H_{n-1}(F) \longrightarrow H_{n-1}(F')\longrightarrow H_{n-1}(F,F')$$
with $H_n(F,F')=0$, $ H_{n-1}(F)$ and $H_{n-1}(F')$ free of rank $\nu$,
and $H_{n-1}(F,F')$ torsion-free.
Then the inclusion of $F'$ in $F$ is an homotopy equivalence.

The second part is based on the $h$-cobordism theorem.
Let $X = f^{-1}(K) \cap B_R \setminus \mathring B_{R'}$, then as
$f$ has no affine critical points in $X$ (because there is no critical values in $K$)
 and $f$ is transversal
to $f^{-1}(K)\cap  S_{R}$ and  to $f^{-1}(K)\cap  S_{R'}$ then by
Ehresmann theorem $f : X \longrightarrow K$ is a fibration. We denote $F
\setminus \mathring F'$  by $F^*$. We get an isomorphism $H_i(\partial
F') \longrightarrow H_i(F^*)$ for all $i$ because $H_i(F^*,\partial
F')=H_i(F,F')=0$. For $n=2$ it implies that $F^*$ is diffeomorphic to a
product $[0,1] \times \partial F'$. For $n>3$ we will use the
$h$-cobordism  theorem to $F^*$ to prove this. We have $\partial F^* =
\partial F' \cup \partial F$; $\partial F'$ and $\partial F$ are simply
connected: if we look at the function $x\mapsto -\|x-x_0\|$ on $f^{-1}(*)$
for a generic $x_0$, then $F=f^{-1}(*)\cap B_R$ and $F'=f^{-1}(*)\cap
B_{R'}$ are obtained by gluing cells of index more or equal to $n-1$. So
their boundary is simply connected. For a similar reason $F^*$ is simply
connected. As we have isomorphisms $H_i(\partial F') \longrightarrow
H_i(F^*)$ and both spaces are simply connected then by Hurewicz-Whitehead
theorem the inclusion of $\partial F'$ in $F^*$ is an homotopy
equivalence. Now $F^*$,  $\partial F'$, $\partial F$ are simply connected,
the inclusion of $\partial F'$ in $F^*$ is an homotopy equivalence and
$F^*$ has real dimension $2n-2 \geq 6$. So by the $h$-cobordism theorem
\cite{M2} $F^*$ is diffeomorphic to the product $[0,1] \times \partial
F'$. Then the fibration $f : X \longrightarrow K$ is differentiably
equivalent
 to the fibration $f : [0,1]\times (f^{-1}(K)\cap S_{R'}) \longrightarrow K$
so the fibrations  $f :  f^{-1}(K)\cap B_{R'}  \longrightarrow K$
and $f :   f^{-1}(K)\cap  B_{R} \longrightarrow K$  are
differentiably equivalent.
\end{proof}

\section{Family of polynomials}

Let $(f_t)_{t\in[0,1]}$ be a family of polynomials that verify hypotheses of theorem \ref{th:A}.

\begin{lemma}[\cite{HP}]
\label{lem:4}
 There exists $R\gg1$ such that for all $t \in [0,1]$
the affine critical points of $f_t$ are in $\mathring B_R$.
\end{lemma}

\begin{proof}
It is enough to prove it on $[0,\tau]$ with $\tau > 0$. We choose $R\gg1$
such that all the affine critical points of $f_0$ are in $\mathring B_R$.
We denote $$\phi_{t} = \frac{\grad_{f_{t}}}{\| \grad_{f_{t}} \|} : S_R
\longrightarrow S_1.$$ Then $\deg \phi_0 = \mu(0)$. For all $x\in S_R$,
$\grad_{f_0}x \not= 0$, and by continuity there exist $\tau > 0$ such
that for $t\in[0,\tau]$ and all $x\in S_R$, $\grad_{f_t}x \not= 0$. Then
the maps $\phi_t$ are homotopic (the homotopy is $\phi : S_R\times
[0,\tau] \longrightarrow S_1$ with $\phi(x,t) = \phi_t(x)$). And then
$\mu(0) = \deg \phi_0 = \deg \phi_t \leq \mu(t)$. If there exists a
family $x(t)\in \Cc^n$ of affine critical   points of $\phi_t$ such that
$\| x(t) \| \rightarrow +\infty$ as $t\rightarrow 0$, then for a
sufficiently small $t$, $x(t) \notin B_R$ and then $\mu(t) > \deg \phi_t$.
It contradicts the hypothesis $\mu(0) = \mu(t)$.
\end{proof}

\begin{lemma}
\label{lem:5} There exists $r\gg 1$ such that the subset $\big\lbrace
(c,t) \in D_r(0)\times [0,1] \ | \ c\in\B(t) \big\rbrace$ is a braid of
$D_r(0)\times [0,1]$.
\end{lemma}

It enables us to choose $* \in \partial D_r(0)$ which is a regular value
for all $f_t$, $t\in [0,1]$. In other words if we enumerate $\B(0)$ as
$\{c_1(0),\ldots,c_m(0)\}$ then there is  continuous functions $c_i :
[0,1] \longrightarrow D_r(0)$ such that for $i\not= j$, $c_i(t) \not=
c_j(t)$. This enables us to identify $\pi_1(\Cc \setminus \B(0),*)$ and
$\pi_1(\Cc \setminus \B(1),*)$.

\begin{proof}
Let $\tau$ be in $[0,1]$ and $c(\tau)$ be a critical value of $f_{\tau}$ then
for all $t$ near $\tau$ there exists a critical value $c(t)$ of $f_t$.
It is an hypothesis for the critical values at infinity
and this fact is well-known for affine critical values as the coefficients of $f_t$
are smooth functions of $t$, see for example \cite[Prop.~2.1]{Br}.

Moreover by the former lemma there can not exist critical values that escape at infinity
\emph{i.e.} a $\tau\in[0,1]$
such that $|c(t)| \rightarrow +\infty$ as $t\rightarrow \tau$.
For affine critical values it is a consequence of the former lemma
(or we can make the same proof as we now will perform for the critical values at
infinity).
For $\Binf(t)$ let us suppose that there is critical values that
escape at infinity. By continuity of the critical values at infinity
with respect to $t$ we can suppose that there is a continuous function
$c_0(t)$ on $]0,\tau]$ ($\tau>0$) with $c_0(t) \in \Binf(t)$ and $|c(t)|\rightarrow +\infty$
as $t\rightarrow 0$.
By continuity of the critical values at infinity, if  $\Binf(0) = \{c_1(0),\ldots,c_p(0)\}$
there exist continuous functions $c_i(t)$ on $[0,\tau]$ such that $c_i(t) \in \Binf(t)$
for all $i= 1,\ldots,p$.
And for a sufficiently small $t>0$, $c_0(t)\not= c_i(t)$ ($i=1,\ldots,p$)
then $\#\Binf(0) < \#\Binf(t)$ which contradicts the constancy of $\#\Binf(t)$.

Finally there can not exist ramification points: suppose that there is a $\tau$
such that $c_i(\tau) = c_j(\tau)$ (and $c_i(t), c_j(t)$ are not equal in a neighborhood of $\tau$).
Then if $c_i(\tau) \in \Baff(\tau) \setminus \Binf(\tau)$
(\emph{resp.} $\Binf(\tau) \setminus \Baff(\tau)$,  $\Binf(\tau) \cap \Baff(\tau)$)
there is jump in $\# \Baff(t)$ (\emph{resp.}  $\# \Binf(t)$, $\# \B(t)$) near $\tau$ which
is impossible  by assumption.
\end{proof}

\bigskip

Let $R_0, K, D_r(0), D_\epsilon(c)$ be the objects of the former section
for the polynomial $f=f_0$. Moreover we suppose that $R_0$ is greater than the $R$
obtained in lemma \ref{lem:4}.

\begin{lemma}
\label{lem:6} There exists $\tau \in ]0,1]$ such that for all $t\in
[0,\tau]$ we have the properties:
\begin{itemize}
  \item $c_i(t) \in D_{\epsilon}(c_i(0))$, $i=1,\ldots,m$;
  \item for all $s \in K$,  $f_t^{-1}(s)$ intersects $S_{R_0}$ transversally.
\end{itemize}
\end{lemma}

\begin{proof}
The first point is just the continuity of the critical values $c_i(t)$.
The second point is the continuity of transversality:
if the property is false then there exists  sequences $t_k \rightarrow 0$,
$x_k \in S_{R_0}$ and $\lambda_k \in \Cc$ such that $\grad_{f_{t_k}}x_k=\lambda_k x_k$.
We can suppose that $(x_k)$ converges (after extraction of a sub-sequence, if necessary).
Then $x_k \rightarrow x \in S_{R_0}$, $\grad_{f_{t_k}}x_k \rightarrow \grad_{f_{0}}x$,
and $\lambda_k = \langle \grad_{f_{t_k}}x_k \, | \, x_k \rangle / \| x_k \|^2 =
\langle \grad_{f_{t_k}}x_k \ | \ x_k \rangle / {R_0}^2$ converges toward $\lambda \in \Cc$.
Then $\grad_{f_{0}}x=\lambda x$ and the intersection is non-transversal.
\end{proof}

\begin{lemma}
\label{lem:7}
The fibrations $f_0 : f_0^{-1}(K) \cap B_{R_0} \longrightarrow K$ and
$f_\tau : f_\tau^{-1}(K) \cap B_{R_0} \longrightarrow K$ are differentiably isomorphic.
\end{lemma}

\begin{proof}
Let
$$F : \Cc^n\times [0,1] \longrightarrow \Cc \times [0,1], \qquad (x,t)\mapsto (f_t(x),t).$$
We want to prove that the fibrations
$F_0 : \Sigma_0 = F^{-1}(K\times\{0\}) \cap (B_{R_0}\times\{0\}) \longrightarrow K$,
$(x,0) \mapsto f_0(x)$
and
$F_\tau : \Sigma_\tau = F^{-1}(K\times\{\tau\}) \cap (B_{R_\tau}\times\{\tau\})
\longrightarrow K$,
$(x,\tau) \mapsto f_\tau(x)$ are differentiably isomorphic.
Let denote $[0,\tau]$ by $I$. Then $F$ has maximal rank on
$F^{-1}(K\times I) \cap (\mathring B_{R_0}\times I)$ and
on the boundary $F^{-1}(K\times I) \cap (S_{R_0}\times I)$.
By Ehresmann theorem $F : F^{-1}(K\times I) \cap (B_{R_0}\times I)
\longrightarrow K\times I$ is a fibration. But we can not argue as in \cite{LR}
since the restriction of $F$ on the set
$\big\lbrace (x,t) \in S_{R_0} \times I \ | \ f_t(x)\in D_r(0)\big\rbrace$
is not a trivial fibration.

As in \cite{HP} we build a vector field that give us a diffeomorphism between the two fibrations
$F_0$ and $F_\tau$.
Let $R_2$ be a real number close to $R_0$ such that $R_2 < R_0$.
On the set $F^{-1}(K\times I) \cap (\cup_{R_2<R<R_0}S_{R}\times I)$
we build a vector field $v_1$ such that for $z\in S_R\times I$ ($R_2<R<R_0$),
$v_1(z)$ is tangent to $S_R \times I$ and we have
$d_z F.v_1(z) = (0,1)$.
On the set $F^{-1}(K\times I) \cap (\mathring B_{R_3}\times I)$
with $R_2<R_3<R_0$ we build a second vector field $v_2$ such that
$d_z F.v_2(z) = (0,1)$, this is possible because $F$ is a submersion on this set.

By gluing these vector fields $v_1$ and $v_2$ by a partition of unity and by integrating
the corresponding vector field we obtain integral curves $p_z : \Rr \longrightarrow
F^{-1}(K\times I)\cap B_{R_0} \times I$ for $z \in \Sigma_0$ such that $p_z(0)=z$
and $p_z(\tau) \in \Sigma_\tau$. It induces a diffeomorphism $\Phi : \Sigma_0 \longrightarrow
\Sigma_\tau$ such that $F_0=F_\tau \circ \Phi$; that makes the fibrations isomorphic.
\end{proof}

\bigskip

\begin{proof}[Proof of theorem \ref{th:A}]
It suffices to prove the theorem for an interval $[0,\tau]$ with $\tau>0$.
We choose $\tau$ as in lemma \ref{lem:6}.
By lemma \ref{lem:0}, $f_0 : f^{-1}(\Cc\setminus\B(0)) \longrightarrow \Cc\setminus \B(0)$
and $f_0 : f_0^{-1}(K) \longrightarrow K$ are differentiably isomorphic fibrations.
Then by lemma \ref{lem:2}, the fibration  $f_0 : f_0^{-1}(K) \longrightarrow K$
is differentiably isomorphic to
$f_0 : f_0^{-1}(K) \cap \mathring B_{R_0} \longrightarrow K$
which is, by lemma \ref{lem:7}  differentiably isomorphic to
$f_\tau : f_\tau^{-1}(K) \cap \mathring B_{R_0} \longrightarrow K$.

By continuity of transversality (lemma \ref{lem:6}) $f_\tau^{-1}(K)$
has transversal intersection with $S_{R_0}$, we choose a large real number $R$
(by lemma \ref{lem:1} applied to $f_\tau$) such that $f_\tau^{-1}(K)$
intersects  $S_{R}$ transversally.
The last fibration is fiber homotopy equivalent to
$f_\tau : f_\tau^{-1}(K) \cap \mathring B_{R} \longrightarrow K$:
 it is the first part of lemma \ref{lem:3} because the fiber
$f_\tau^{-1}(*)\cap \mathring B_{R_0}$ is homotopic to a wedge
of $\mu(0)+\lambda(0)$ circles and the fiber $f_\tau^{-1}(*)\cap \mathring B_{R}$
is homotopic to a wedge of $\mu(\tau)+\lambda(\tau)$ circles; as
$\mu(0)+\lambda(0)=\mu(\tau)+\lambda(\tau)$ we get the desired conclusion.
Moreover for $n\not=3$ by the second part of lemma \ref{lem:3} the fibrations are
differentiably isomorphic.

Applying lemma \ref{lem:2} and \ref{lem:0} to $f_\tau$ this fibration is differentiably
isomorphic
to $f_\tau : f_\tau^{-1}(\Cc\setminus\B(\tau)) \longrightarrow \Cc\setminus\B(\tau)$.
As a conclusion the fibrations $f_0 : f_0^{-1}(\Cc\setminus\B(0))
\longrightarrow \Cc\setminus\B(0)$ and
$f_\tau : f_\tau^{-1}(\Cc\setminus\B(\tau)) \longrightarrow \Cc\setminus\B(\tau) $ are
fiber homotopy equivalent,
and for $n\not=3$ are differentiably isomorphic
\end{proof}

\section{Around affine singularities}

We now work with $t\in [0,1]$. We suppose  in this paragraph that the
critical values $\B(t)$ depend analytically on $t\in [0,1]$. This enables
us to construct a diffeomorphism: $$\chi : \Cc \times [0,1]
\longrightarrow \Cc\times [0,1], \quad \text{ with } \chi(x,t) =
(\chi_t(x),t),$$ such that $\chi_0 = \id$ and $\chi_t(\B(t))= \B(0)$. We
denote $\chi_1$ by $\Psi$, so that $\Psi : \Cc \longrightarrow \Cc$
verify $\Psi(\B(1)) = \B(0)$. Moreover we can suppose that $\chi_t$ is
equal to $\id$ on $\Cc \setminus D_r(0)$ this is possible because for all
$t\in[0,1]$, $\B(t)\subset D_r(0)$. Finally $\chi$ defines a vector field
$w$ of $\Cc \times [0,1]$ by $\frac{\partial \chi}{\partial t}$.

\bigskip

We need a non-splitting of the affine singularity, this principle has been
proved by C.~Has~Bey (\cite{HB}, $n=2$) and by F.~Lazzeri (\cite{La}, for all $n$).

\begin{lemma}
\label{lem:8} Let $x(\tau)$ be an affine singular point of $f_\tau$ and
let $U_\tau$ be an open neighborhood of $x(\tau)$ in $\Cc^n$ such that
$x(\tau)$ is the only affine singular point of $f_\tau$ in $U_\tau$.
Suppose that for all $t$ closed to $\tau$, the restriction of $f_t$ to
$U_\tau$ has only one critical value.
 Then for all $t$ sufficiently closed to $\tau$, there is
one, and only one, affine singular point of $f_t$ contained in $U_\tau$.
\end{lemma}

So we can enumerate the singularities:
 if we denote the affine singularities of
$f_0$ by $\{ x_i(0) \}_{i\in J}$ then there is continuous functions $x_i
: [0,1] \longrightarrow \Cc^n$ such that $\{x_i(t)\}_{i\in J}$ is the set
of affine singularities of $f_t$. Let us notice that there can be two
distinct singular points of $f_t$ with the same critical value. We
suppose that $(f_t)$ verifies the hypotheses of theorem \ref{th:A}, that
$n\not=3$, and $\B(t)$ depends analytically on $t$. This and the former
lemma imply that for all $t\in[0,1]$ the local Milnor number of $f_t$ at
$x(t)$ is equal to the local Milnor number of $f_0$ at $x(0)$. The
improved version of L\^e-Ramanujam theorem by J.G.~Timourian \cite{Tim}
for a family of germs with  constant local Milnor numbers proves that
$(f_t)$ is locally a product family.

\begin{theorem}[L\^e-Ramanujam-Timourian]
Let $x(t)$ be a singular points of $f_t$. There exists $U_{t}$, $V_{t}$
neighborhoods of $x(t)$, $f_t(x(t))$ respectively and an homeomorphism
$\Omega^\ine$ such that if $U = \bigcup_{t\in[0,1]}U_t\times\{t\}$ and $V
= \bigcup_{t\in[0,1]}V_t\times\{t\}$ the following diagram commutes: $$
\xymatrix{
{}U \ar[r]^-{\Omega^\ine} \ar[d]_-{F}  & U_0\ar[d]^-{f_0\times\id} \times [0,1] \\
  V  \ar[r]_-\chi       & V_0 \times [0,1].
}
$$
\end{theorem}

In particular it proves that the polynomials $f_0$ and $f_1$
are locally topologically equivalent: we get an homeomorphism $\Phi_\ine$
 such that the following diagram commutes:
$$
\xymatrix{
{}U_1\ar[r]^-{\Phi_\ine}\ar[d]_-{f_1}  & U_0\ar[d]^-{f_0} \\
V_1 \ar[r]_-\Psi        & V_0. } $$ By lemma \ref{lem:4} we know that for
all $t \in [0,1]$, $\B(t) \subset D_r(0)$. We extend the definition of
$R_0$ and $R_1$ to all $f_t$. Be continuity of transversality and
compactness of $[0,1]$ we choose $R_1$ such that
 $$\forall c\in\B(0)\quad \forall R\ge R_1 \quad f_0^{-1}(c) \pitchfork S_{R}
 \quad \text{ and } \quad \forall t\in[0,1]\quad \forall c\in\B(t) \quad  f_t^{-1}(c) \pitchfork S_{R_1}.$$
 For a sufficiently small $\epsilon$ we denote
$$K(0) = D_r(0) \setminus \bigcup_{c\in\Binf(0)} D_\epsilon(c),
\qquad K(t) = \chi_t^{-1}(K(0))$$ and we choose $R_0\ge R_1$ such
that
  $$\forall s\in K(t) \quad \forall   R \ge R_0 \quad f_0^{-1}(s) \pitchfork S_{R}
  \quad \text{ and } \quad
  \forall t\in[0,1] \quad \forall s\in K(t)  \quad f_t^{-1}(s) \pitchfork S_{R_0}.$$
We denote $$B'_t=
B_{R_1} \cup \big( f_t^{-1}(K(t))  \cap \mathring B_{R_0}\big),
\quad t\in [0,1].$$

\begin{lemma}
\label{lem:9}
There exists an homeomorphism $\Phi$
such that we have the commutative diagram:
$$
\xymatrix{
{}B'_1\ar[r]^-{\Phi}\ar[d]_-{f_1}  & B'_0\ar[d]^-{f_0} \\
D_r(0) \ar[r]_-\Psi        & D_r(0).
}
$$
\end{lemma}

\begin{proof}
We denote by $U'_t$ a  neighborhood of $x(t)$ such that  $\bar {U'_t}
\subset U_t$. We denote by $\mathcal{U}_t$ (resp. $\mathcal{U}'_t$), the
union (on the affine singular points of $f_t$ ) of the $U_t$ (resp.
$U_t'$). We set $$B''_t = B'_t \setminus \mathcal{U}'_t, \quad t\in
[0,1].$$ We can extend the homeomorphism $\Phi$ of lemma \ref{lem:7} to
$\Phi_\out: B''_1 \longrightarrow  B''_0$. We just have to extend the
vector field of lemma \ref{lem:7} to a new vector field denoted by $v'$
such that
\begin{itemize}
  \item $v'$ is tangent to $\partial \mathcal{U}'_t$,
  \item $v'$ is tangent to $S_{R_1}\times [0,1]$ on $F^{-1}(D_r(0) \setminus K(t) \times \{ t \})$
for all $t\in [0,1]$,
  \item $v'$ is tangent to $S_{R_0}\times [0,1]$ on $F^{-1}(K(t)\times \{ t \})$
for all $t$.
  \item $d_zF.v'(z) = w(F(z))$ for all $z \in \bigcup_{t\in[0,1]}B''_t\times\{t\}$, which means that $\Phi_\out$
respect the fibrations.
\end{itemize}
If we set $B'' = \bigcup_{t\in[0,1]}B''_t\times\{t\}$ the integration of $v'$
gives $\Omega_\out$ and $\Phi_\out$ such that:
$$
\xymatrix{
{} B'' \ar[r]^-{\Omega^\out} \ar[d]_-{F}  & B''_0\ar[d]^-{f_0\times\id} \times [0,1] \\
  \Cc\times[0,1]  \ar[r]_-\chi       & \Cc \times [0,1],
}
\qquad \qquad
\xymatrix{
{}B''_1\ar[r]^-{\Phi_\out}\ar[d]_-{f_1}  & B''_0\ar[d]^-{f_0} \\
D_r(0) \ar[r]_-\Psi        & D_r(0).
}
$$

We now explain how to glue $\Phi_\ine$ and $\Phi_\out$ together.
We can suppose that there exists spheres $S_t$ centered at the
singularities $x(t)$ such that if $S = \bigcup_{t\in[0,1]}S_t
\times\{t\}$ $\Omega^\ine : S \longrightarrow S_0 \times [0,1]$
and $\Omega^\out : S \longrightarrow S_0 \times [0,1]$. It defines
$\Omega^\ine_t : S_t \longrightarrow S_0$ and $\Omega^\out_t : S_t
\longrightarrow S_0$. Now we define
$$\Theta_t : S_1 \longrightarrow S_0, \quad \Theta_t
   = \Omega^\ine_t \circ (\Omega^\out_t)^{-1} \circ \Phi_\out.$$
Then $\Theta_0 = \Phi_\out$ and $\Theta_1 = \Phi_\ine$. On a set
homeomophic to $S\times [0,1]$ included in $\bigcup_{t\in[0,1]}U_t
\setminus U'_t$ we glue $\Phi_\ine$ to $\Phi_\out$, moreover this gluing
respect the fibrations $f_0$ and $f_1$. We end by doing this construction
for all affine singular points.
\end{proof}

\begin{proof}[Proof of theorem \ref{th:B}]
We firstly prove that affine critical values are analytic functions of
$t$. Let $c(0) \in \Baff(0)$, the set $\big\{ (c(t), t) \ | \ t\in[0,1]
\big\}$ is a real algebraic subset of $\Cc \times [0,1]$ as all affine
critical points are contained in $B_{R_0}$ (lemma \ref{lem:4}). In fact
there is a polynomial $P \in \Cc[x,t]$ such that $(c=0)$ is equal to
$(P=0) \cap \Cc \times [0,1]$. Because the set of critical values is a
braid of $\Cc \times [0,1]$ (lemma \ref{lem:5}) then $c : [0,1]
\longrightarrow \Cc$ is a smooth analytic function.

If we suppose that $\Binf(t) = \varnothing$ for all $t\in[0,1]$
then by lemma \ref{lem:2} we can extend $\Phi : B'_1 \longrightarrow B'_0$
to $\Phi : f_1^{-1}(D_r(0)) \longrightarrow f_0^{-1}(D_r(0))$.
And as $\B(t) \subset D_r(0)$ by a lemma similar to lemma \ref{lem:0}
we can extend the homeomorphism to the whole space.
\end{proof}

\section{Polynomials in two variables}

We set $n=2$. Let $f_t : \Cc^2 \longrightarrow \Cc$ such that the
coefficient of this family are algebraic in $t$. We suppose that the
integers $\mu(t)$, $\lambda(t)$, $\#\B(t)$, $\#\Baff(t)$, $\#\Binf(t)$ do
not depend on $t\in[0,1]$. We also suppose the $\deg f_t$ does not depend
on $t$.

We recall a result of L.~Fourrier \cite{Fo}. Let  $f:\Cc^2
\longrightarrow \Cc$ with set of critical values at infinity $\B_\infty$.
Let $* \notin \B$ and $Z = f^{-1}(*) \cup \bigcup_{c\in\Binf}f^{-1}(c)$.
The \defi{total link of $f$} is $L_f = Z \cap S_R$ for a sufficiently
large $R$. To $f$ we associate a resolution $\phi : \Sigma
\longrightarrow \Pp^1$, the components of the divisor of this resolution
on which $\phi$ is surjective are the \defi{dicritical components}. For
each dicritical component $D$ we have a branched covering $\phi : D
\longrightarrow \Pp^1$. If the set of dicritical components is $D_\dic$
we then have the restriction of $\phi$, $\phi_\dic : D_\dic
\longrightarrow \Pp^1$. The \defi{$0$-monodromy representation} is the
representation $$\pi_1(\Cc\setminus \B) \longrightarrow \Aut
\big(\phi_\dic^{-1}(*)\big).$$
\begin{theorem}[Fourrier]
Let $f,g$ be complex polynomials in two variables with equivalent
$0$-monodromy representations and equivalent total links then there exist
homeomorphisms $\Phi_\infty$ and $\Psi_\infty$ and compact sets $C,C'$ of
$\Cc^2$ that make the diagram commuting: $$ \xymatrix{
{}\Cc^2\setminus C \ar[r]^-{\Phi_\infty}\ar[d]_-{f}  & \Cc^2\setminus C' \ar[d]^-{g} \\
  \Cc \ar[r]_-{\Psi_\infty}        & \Cc.
}
$$
\end{theorem}

For our family $(f_t)$, by  theorem \ref{th:A} we know that the
geometric monodromy representations are all equivalent, it implies
that all the $0$-monodromy representations of $(f_t)$ are
equivalent. Moreover if we suppose that for any $t,t'\in[0,1]$ the
total links $L_{f_t}$ and $L_{f_{t'}}$ are equivalent, then by the
former theorem the polynomials $f_t$ and $f_{t'}$ are
topologically equivalent out of some compact sets of $\Cc^2$. We
need a result a bit stronger which can be proved by similar
arguments than in \cite{Fo} and we will omit the proof:
\begin{lemma}
\label{lem:Fo}
Let $(f_t)_{t\in[0,1]}$ be a polynomial family such
that the coefficients are algebraic functions of $t$. We suppose
that the $0$-monodromy representations and the total links are all
equivalent. Then there exists compact sets $C(t)$ of $\Cc^2$ and
an homeomorphism $\Omega^\infty$ such that if $\mathcal{C} =
\bigcup_{t\in[0,1]} C(t)\times\{t\}$ we have a commutative
diagram: $$ \xymatrix{ {} \Cc^2\times [0,1] \setminus \mathcal{C}
\ar[r]^-{\Omega^\infty} \ar[d]_-{F}  &
             \Cc^2 \setminus C(0) \times [0,1] \ar[d]^-{f_0\times\id} \\
  \Cc\times[0,1]  \ar[r]_-\chi       & \Cc \times [0,1].
}
$$
\end{lemma}

We now prove a strong version of the continuity of critical
values.
\begin{lemma}
\label{lem:mucst}
 The critical values are smooth analytic functions of
$t$. Moreover for $c(t) \in \B(t)$, the integer $\mu_{c(t)}$ and
$\lambda_{c(t)}$ do not depend on $t\in[0,1]$.
\end{lemma}

\begin{proof}
For affine critical values, refer to the proof of theorem
\ref{th:B}. The constancy of $\mu_{c(t)}$ is a consequence of
lemma \ref{lem:4} and lemma \ref{lem:8}. For critical values at
infinity we need a result of \cite{H} and \cite{HP} that enables
to calculate critical values and Milnor numbers at infinity. As
$\deg f_t$ is constant we can suppose that this degree is $\deg_y
f_t$. Let denote $\Delta (x,s,t)$ the discriminant $\Disc_y
(f_t(x,y)-s)$ with respect to $y$. We write
  $$\Delta(x,s,t) = q_1(s,t)x^{k(t)}+q_2(s,t)x^{k(t)-1}+\cdots $$
First of all $\Delta$ has constant degree $k(t)$ in $x$ because
$k(t)= \mu(t) + \lambda(t) + \deg f_t -1$ (see \cite{HP}).
Secondly by \cite{H} we have
 $$\Binf(t) = \big\lbrace s \ |\ q_1(s,t) = 0 \big\rbrace$$
then we see that critical values at infinity depend continuously
on $t$ and that critical values at infinity are a real algebraic
subset of $\Cc \times [0,1]$. For the analicity we end as in the
proof of theorem \ref{th:B}. Finally, for a fixed $t$, we have
that $\lambda_{c} = k(t) - \deg_x \Delta(x,c,t)$. In other words
$q_i(t,c)$ is zero for $i = 1,\ldots,\lambda_c$ and non-zero for
$i = \lambda_c + 1$. For $c(t)\in \Binf(t)$ we now prove that
$\lambda_{c(t)}$ is constant. The former formula proves that
$\lambda_{c(t)}$ is constant except for all but finitely many
$\tau\in[0,1]$ for which $\lambda_{c(\tau)} \ge \lambda_{c(t)}$.
But if $\lambda_{c(\tau)}
> \lambda_{c(t)}$ then $\lambda(\tau) = \sum_{c\in\Binf(\tau)}
\lambda_{c} > \sum_{c\in\Binf(t)} \lambda_{c} = \lambda(t)$ which
contradicts the hypotheses.
\end{proof}

To apply lemma \ref{lem:Fo} we need to prove:
\begin{lemma}
For any $t,t'\in[0,1]$ the total links $L_{f_t}$ and $L_{f_t'}$
are equivalent.
\end{lemma}

\begin{proof}
The problem is similar to the one of \cite{LR} and to lemma
\ref{lem:3}. For a value $c(t)$ in $\Binf(t)$ or equal to $*$, we
have that the link at infinity $f_0^{-1}(c(0))\cap S_{R_1}$ is
equivalent to the link $f_1^{-1}(c(1))\cap S_{R_1}$ (lemma
\ref{lem:9}). But $f_1^{-1}(c(1))\cap S_{R_1}$ is not necessarily
the link at infinity for $f_1^{-1}(c(1))$. We now prove this fact;
let denote $c=c(1)$. Let $R_2\ge R_1$ such that for all $R\ge
R_2$, $f_1^{-1}(c)\pitchfork S_{R}$, then $f_1^{-1}(c)\cap
S_{R_2}$ is the link at infinity of $f_1^{-1}(c)$. We choose
$\eta$, $0<\eta\ll 1$ such that $f_1^{-1}(D_\eta(c))$ has
transversal intersection with $S_{R_1}$ and $S_{R_2}$ and such
that $f_1^{-1}(\partial D_\eta(c))$ has transversal intersection
with all $S_R$, $R\in [R_1,R_2]$. Notice that $\eta$ is much
smaller than the $\epsilon$ of the former paragraphs and that
$f_1^{-1}(s)\cap S_{R_2}$ is \textbf{not} the link at infinity of
$f_1^{-1}(s)$ for $s\in \partial D_\eta(c)$. We fix $R_0$ smaller
than $R_1$ such that $f_1^{-1}(D_\eta(c))$ has transversal
intersection with $S_{R_0}$. We denote $f_1^{-1}(D_\eta(c)) \cap
B_{R_i} \setminus \mathring B_{R_0}$ by $\mathcal{A}_i$, $i =
1,2$. The proof is now similar to the one of lemma \ref{lem:3}.
Let $A_1$ and $A_2$ be connected components of $\mathcal{A}_1$ and
$\mathcal{A}_2$ with $A_1 \subset A_2$. By Ehresmann theorem, we
have fibrations $f_1: A_1 \longrightarrow D_\eta(c)$, $f_1: A_2
\longrightarrow D_\eta(c)$. From one hand $f_1^{-1}(c)\cap
B_{R_1}$ has the homotopy type of a wedge of
$\mu+\lambda-\mu_{c(0)}-\lambda_{c(0)}$ circles, because
$f_1^{-1}(c)\cap B_{R_1}$ is diffeomorphic to $f_1^{-1}(c(0))\cap
B_{R_1}$ with Euler characteristic
$1-\mu-\lambda+\mu_{c(0)}+\lambda_{c(0)}$ by Suzuki formula. From
the other hand $f_1^{-1}(c)\cap B_{R_2}$ has the homotopy type of
a wedge of $\mu+\lambda-\mu_{c(1)}-\lambda_{c(1)}$ circles by
Suzuki formula. By lemma \ref{lem:mucst} we have that
$\mu_{c(0)}+\lambda_{c(0)}= \mu_{c(1)}+\lambda_{c(1)}$, with
$c=c(1)$, so the fiber $f_1^{-1}(c)\cap B_{R_1}$ and
$f_1^{-1}(c)\cap B_{R_2}$ are homotopic, it implies that the
fibrations $f_1: A_1 \longrightarrow D_\eta(c)$ and  $f_1: A_2
\longrightarrow D_\eta(c)$ are  fiber homotopy equivalent, and
even more are diffeomorphic. It provides a diffeomorphism $\Xi :
A_1\cap S_{R_1}=A_2\cap S_{R_1} \longrightarrow A_2\cap S_{R_2}$
and we can suppose that $\Xi(f_1^{-1}(c)\cap A_1 \cap S_{R_1})$ is
equal to $f_1^{-1}(c)\cap A_1 \cap S_{R_1}$. By doing this for all
connected components of $\mathcal{A}_1$, $\mathcal{A}_2$, for all
values $c\in \Binf(1)\cup\{*\}$ and by extending $\Xi$ to the
whole spheres we get a diffeomorphism $\Xi : S_{R_1}
\longrightarrow S_{R_2}$ such that $\Xi(f_1^{-1}(c)\cap S_{R_1})=
f_1^{-1}(c)\cap S_{R_2}$ for all $c\in \Binf(1)\cup\{*\}$. Then
the total link for $f_0$ and $f_1$ are equivalent.
\end{proof}

\begin{proof}[Proof of theorem \ref{th:C}]
 By lemma \ref{lem:Fo} we have a trivialization $\Omega^\infty
: \Cc^2\times [0,1] \setminus \mathcal{C}
  \longrightarrow  \Cc^2 \setminus C(0)\times [0,1]$.
We can choose the $R_1$ (before lemma \ref{lem:9}) such that $\mathring
C(t) \subset B_{R_1}$. And then the proof of this lemma gives us an
$\Omega^\out : \bigcup_{t\in[0,1]} B''(t)\times\{t\} \longrightarrow
B''(0)\times [0,1]$. By gluing $\Omega^\out$ and $\Omega^\infty$ as in
this proof we obtain $\Phi : \Cc^2 \longrightarrow \Cc^2$ such that: $$
\xymatrix{
{}\Cc^2 \ar[r]^-\Phi\ar[d]_-{f_1}  & \Cc^2 \ar[d]^-{f_0} \\
  \Cc \ar[r]_-\Psi        & \Cc.
}
$$
Then $f_0$ and $f_1$ are topologically equivalent.
\end{proof}

\section{Continuity of the critical values at infinity}

\begin{lemma}
\label{prop:D} Let $(f_t)_{t\in[0,1]}$ be a family of polynomials such
that the coefficients are polynomials in $t$. We suppose that:
\begin{itemize}
\item the total affine Milnor number $\mu(t)$ is constant;
\item the degree $\deg f_t$ is a constant;
\item the set of critical points at infinity $\Sigma(t)$ is finite
and does not vary: $\Sigma(t)=\Sigma$;
\item for all $(x:0)\in \Sigma$, the generic Milnor number
$\nu_{\check x}(t)$ is independent of $t$.
\end{itemize}
Then the critical values at infinity depend continuously on $t$,
\emph{i.e.} if $c(t_0) \in \Binf(t_0)$ then for all $t$ near $t_0$ there
exists $c(t)$ near $c(t_0)$ such that $c(t) \in \Binf(t)$.
\end{lemma}

Let $f$ be a polynomial. For $x \in \Cc^n$ we have $(x:1)$ in $\Pp^n$ and
if $x_n \not=0$ we divide $x$ by $x_n$ to obtain local coordinates at
infinity $(\check{x}',x_0)$. The following lemma explains the link between
the critical points of $f$ and those of $F_c$. It uses Euler relation for
the homogeneous polynomial of $f$ of degree $d$.
\begin{lemma}\
\label{lem:13}
\begin{itemize}
\item $F_c$ has a critical point $(\check{x}',x_0)$ with $x_0 \not= 0$ of critical value $0$
if and only if $f$ has a critical point $x$ with critical value $c$.
\item $F_c$ has a critical point $(\check{x}',0)$ of critical value $0$
if and only if  $(x:0)\in \Sigma$.
\end{itemize}
\end{lemma}

\begin{proof}[Proof of lemma \ref{prop:D}]
We suppose that critical values at infinity are \emph{not} continuous
functions of $t$. Then there exists $(t_0,c_0)$ such that $c_0 \in
\Binf(t_0)$ and for all $(t,c)$ in a neighborhood of  $(t_0,c_0)$, we
have $c \notin \Binf(t)$. Let $P$ be the point of irregularity at
infinity for  $(t_0,c_0)$. Then  $\mu_P(F_{t_0,c_0}) > \mu_P(F_{t_0,c})$
($c\not= c_0$) by definition of $c_0 \in \Binf(t_0)$ and by
semi-continuity of the local Milnor number at $P$ we have
$\nu_P(t_0)=\mu_P(F_{t_0,c}) \ge \mu_P(F_{t,c})=\nu_P(t)$, $(t,c) \not=
(t_0,c_0)$.

We consider $t$ as a complex parameter.
By continuity of the critical points and by conservation of the Milnor number
for $(t,c) \not= (t_0,c_0)$ we have critical points $M(t,c)$ near $P$ of $F_{t,c}$
that are not equal to $P$.
This fact uses that $\deg f_t$ is a constant, in order to prove that
$F_{t,c}$ depends continuously of $t$.

Let denote by $V'$ the algebraic variety of $\Cc^3\times \Cc^n$ defined
by $(t,c,s,x) \in V'$ if and only if $F_{t,c}$ has a critical point $x$
with critical value $s$ (the equations are $\grad F_{t,c}(x) = 0,
F_{t,c}(x) = s$). If $\mu_P(F_{t,c}) > 0$ for a generic $(t,c)$ then
$\big\{ (t,c,0,P) \ | \ (t,c)\in \Cc^2 \big\} $ is a subvariety of $V'$.
We define $V$ to be the closure of $V'$ minus this subvariety. Then for a
generic  $(t,c)$, $(t,c,0,P) \notin V$. We call $\pi : \Cc^3\times \Cc^n
\longrightarrow \Cc^3$ the projection on the first factor. We set $W =
\pi(V)$. Then $W$ is locally an algebraic variety around $(t_0,c_0,0)$.
For each $(t,c)$ there is a non-zero finite number of values $s$ such that
$(t,c,s) \in W$. So $W$ is locally an equi-dimensional variety of
codimension $1$. Then it is a germ of hypersurface of $\Cc^3$. Let
$P(t,c,s)$ be the polynomial that defines $W$ locally. We set $Q(t,c) =
P(t,c,0)$. As $Q(t_0,c_0) = 0$ then in all neighborhoods of $(t_0,c_0)$
there exists $(t,c) \not= (t_0,c_0)$ such that $Q(t,c)=0$. Moreover there
are solutions for $t$ a real number near $t_0$.

Then for $(t,c)\not= (t_0,c_0)$ we have that: $Q(t,c)=0$ if and only if
$F_{t,c}$ has a critical point $M(t,c) \not= P$ with critical value $0$.
The point $M(t,c)$ is not equal to $P$ because as $t\not= t_0$,
$(t,c,0,P) \notin V$: it uses that $c\notin \Binf(t)$ for $t\not=t_0$, and
that $\nu_{P}(t)=\nu_{P}(t_0)$. Let us notice that $M(t,c) \rightarrow P$
as $(t,c) \rightarrow (t_0,c_0)$.

We end the proof be studying the different cases:
\begin{itemize}
\item if we have $M(t,c)$ in $\mathcal{H}_\infty$ (of equation $(x_0=0)$) then $M(t,c) \in \Sigma$
which provides a contradiction because then it is equal to $P$;
\item if we have points $M(t,c)$, not in $\mathcal{H}_\infty$, with $t\not= t_0$ then
there are  affine critical points $M'(t,c)$ of $f_t$ (lemma
\ref{lem:13}),  and as  $M(t,c)$ tends towards $P$ (as $(t,c)$ tends
towards $(t_0,c_0)$) we have that $M'(t,c)$ escapes at infinity, it
contradicts the fact that critical points of $f_t$ are bounded (lemma
\ref{lem:4}).
\item if we have points $M(t_0,c)$, not in $\mathcal{H}_\infty$, then
there is  infinitely many affine critical points for $f_{t_0}$, which is impossible
since the singularities of $f_{t_0}$ are isolated.
\end{itemize}
\end{proof}

\section{Examples}

\begin{example}
Let $f_t = x(x^2y+tx+1)$. Then $\Baff(t) = \varnothing$, $\Binf(t)=\{0\}$,
$\lambda(t)=1$ and $\deg f_t = 4$. The by theorem \ref{th:C}, $f_0$ and
$f_1$ are topologically equivalent. These are examples of polynomials
that are topologically but not algebraically equivalent, see \cite{Bo}.
\end{example}

\begin{example}
Let $f_t=(x+t)(xy+1)$. Then $f_0$ and $f_1$ are not topologically
equivalent. One has $\Binf(t) = \varnothing$, $\Baff(t) = \{0,t\}$ for
$t\not=0$, but $\Binf(0)=\{0\}$, $\Baff(0)=\varnothing$. In fact the two
affine critical points for $f_t$ ``escape at infinity'' as $t$ tends
towards $0$.
\end{example}

\begin{example}
Let $f_t = x\big( x(y+tx^2)+1\big)$. Then $f_0$ is topologically equivalent
to $f_1$. We have for all $t\in[0,1]$, $\Baff(t) = \varnothing$, $\Binf(t)=\{0\}$,
and $\lambda(t)=1$, but $\deg f_t = 4$ for $t\not=0$ while $\deg f_0=3$.
\end{example}


{

}

\bigskip
\noindent
\textsf{Arnaud Bodin}\\
\textsf{Centre de Recerca Matem\`atica,  Apartat 50, 08193 Bellatera, Spain}\\
 \texttt{abodin@crm.es}

\end{document}